%% file: relu_quadrature.tex
\documentclass[11pt,reqno]{amsart}
\usepackage[T1]{fontenc}
\usepackage{lmodern}
\usepackage{amsmath,amssymb,amsthm}
\usepackage{booktabs,graphicx}
\usepackage[a4paper,margin=27mm]{geometry}
\usepackage[hidelinks]{hyperref}
\usepackage{cleveref}
\usepackage[
    doi=false,
    isbn=false,
    sortcites,
    backend=biber,
    giveninits=true,
    maxbibnames=99,
    ]{biblatex} 
\usepackage{tikz}
\usetikzlibrary{fit}

\input{mymacros.tex}
\newcommand{\nn}{\mathcal{N}}
\newcommand{\realization}{\mathcal{R}}

\newcommand{\Q}{\mathbb{Q}}

\title[Quadrature by model counting]{Fast quadrature of neural networks via \texorpdfstring{$\#\mathrm{SAT}$}{\#SAT}}
\hypersetup{pdftitle={Quadrature of ReLU networks by model counting}}
\author{Michael Feischl, Benedikt Gr\"a\ss le, Fabian Zehetgruber}
\thanks{The technical details were worked out using AI based on the main ideas/arguments from the authors and we provide a LEAN certificate of the $d=1$ case of Theorem~\ref{thm:counting_quadrature} at \url{https://github.com/michaelFeischl/SatQuad}.
Funded by the Deutsche Forschungsgemeinschaft (DFG, German Research Foundation) -- Project-ID 258734477 -- SFB 1173, the Austrian Science Fund (FWF)
under the special research program Taming complexity in PDE systems (grant SFB F65) as well as projects I6667-N and P36150. This work was also supported by the ERC grant OPTIMAL (grant agreement No. 101125225). }
\begin{document}

\begin{abstract}
       We show that any algorithm for counting satisfying assignments of boolean formulas (model counting, $\#{\rm SAT}$) can be used to approximate integrals where the integrand is the realization of a neural network. The overhead to translate the quadrature problem into a model counting problem is negligible and new subexponential algorithms for counting 3-SAT solutions translate to super algebraic quadrature algorithms for neural networks without curse of dimensionality. 
\end{abstract}
\maketitle

\section{Introduction}
\label{sec:introduction}

Numerical quadrature approximates integrals by weighted sums of function
values. Classical rules exploit smoothness to obtain high accuracy~\cite{quadrature} and suffers from the curse of dimensionality in higher dimensions. Monte Carlo integration provides dimension independent error estimates but converges slowly in general.
Quasi-Monte Carlo methods can improve this rate of convergence under suitable regularity
and weight assumptions~\cite{qmc}, while sparse grids reduce the cost of
tensor-product approximation for functions with mixed
regularity~\cite{BungartzGriebel2004}. The effectiveness of these methods
therefore depends heavily on the smoothness of the integrand as well as its
dimension.

Integrating neural networks is relevant, in particular, to the evaluation
of variational losses in neural-network methods for partial differential
equations. One usually does not want to impose strong smoothness requirements on the neural network class, as this impacts training effectiveness.

Rivera et al.~\cite{RIVERA2022114710} discuss quadrature errors
in this setting and investigate Monte Carlo sampling, adaptive integration,
polynomial approximation, and regularization. Chen et
al.~\cite{chen2021qmc} analyze quasi-Monte Carlo sampling for the deep Ritz
method. Approaches that use the network representation directly include forward integration algorithms with a corrective
procedure~\cite{liu2024neuralnetworksintegrable} and an adaptive method that gives certified error bounds~\cite{certquad}. 

In this work, we exploit the finite arithmetic description of a ReLU
network with rational parameters to evaluate a tensor midpoint sum by
Boolean model counting. We encode the entire sum in one $3$-CNF formula
whose satisfying assignments determine the quadrature value exactly.
A CNF formula is a boolean formula in finitely many variables of the form
\begin{align*}
(\circ x_1\vee \ldots \vee \circ x_{n_1})\wedge (\circ x_{n_1+1}\vee \ldots \vee \circ x_{n_2}) \wedge \ldots,
\end{align*}
where $\circ$ is either the identity or the negation $\lnot$.
Roughly, the idea is to encode the input and the output of the neural network in binary up to a given precision $p$ and rewrite the resulting map from $\{0,1\}^{pd}\to\{0,1\}^p$ as a CNF formula (we can even show that each clause of the formula contains at most three variables).
The construction avoids explicit enumeration of the grid and has size
polynomial in the network description and the logarithm of the grid
resolution. We combine this reduction with a second-order quadrature
error estimate controlled by the network weights. The resulting runtime
bounds depend on the counting algorithm. Together with new subexponential counting algorithms they yield superalgebraic
convergence for each fixed network. Section~\ref{sec:counting_quadrature}
develops the reduction and its analysis, and Section~\ref{sec:rates}
compares the method with uniform Monte Carlo integration.

\section{Deterministic quadrature via \texorpdfstring{$\#\mathrm{SAT}$}{\#SAT}}
\label{sec:counting_quadrature}

A quadrature sum can be evaluated by counting the satisfying assignments
of one Boolean formula, without enumerating its grid points.
The reduction and its error estimate are independent of any particular
counting algorithm.

For a scalar-valued ReLU feed-forward network $\Phi$ with input dimension $d$ and depth $L$ define its realization $\realization_\Phi$ as
\begin{align*}
\realization_\Phi(x) = \bW_L\phi(\bW_{L-1}\phi(\cdots \bW_2\phi(\bW_1x+\bb_1)+\bb_2\cdots)+\bb_{L-1})+\bb_L\in \R,
\end{align*} 
with weight matrices $\bW_\ell$ and bias vectors $\bb_\ell$ of dimensions bounded by $w_\ell\in \N$ for all $\ell=1,\ldots, L$.
With $|\bW_\ell|$ denoting the elementwise absolute value matrix, define
\begin{align*}
    K_\Phi&:=\prod_{\ell=1}^{L}\norm{\bW_\ell}{\infty},
    & (c_1,\ldots,c_d)&:=|\bW_{L}|\cdots|\bW_1|\in \R^{1\times d}.
\end{align*}
Then $\sum_{a=1}^d c_a\leq K_\Phi$ by
submultiplicativity of the induced maximum norm. For $m\in\N$, define the tensor grid and tensor quadrature formula
\begin{align}\label{eq:counting_midpoints}
    J_{m,d}&:=\{0,\ldots,2^m-1\}^d,
    &\bx_{\boldsymbol j}&:=2^{-m}(\boldsymbol j+\boldsymbol{1/2}),
    &Q_m(\Phi)&:=2^{-md}\sum_{\boldsymbol j\in J_{m,d}}
        \realization_\Phi(\bx_{\boldsymbol j}).
\end{align}

\begin{lemma}\label{lem:counting_midpoint_error}
For every scalar-valued ReLU network $\Phi$ and $m\in\N$,
\begin{align}\label{eq:counting_midpoint_error}
    \left|\int_{[0,1]^d}\realization_\Phi(x)\,dx-Q_m(\Phi)\right|
    \leq\frac{2^{-2m}}4\sum_{a=1}^dc_a
    \leq\frac{K_\Phi}{4}2^{-2m}.
\end{align}
\end{lemma}
\begin{proof}
For a continuous piecewise affine function $g:[0,1]\to\R$, write
\[
    V_*(g):=\operatorname{Var}_{(0,1)}(g')
                +|g'(0+)|+|g'(1-)|.
\]
Here the variation of the piecewise constant derivative is the sum of
the absolute jumps at its interior breakpoints. Equivalently, $V_*(g)$
is the variation of the derivative of the constant extension of $g$
to $\R$. This is the intrinsic variability used
in~\cite[Definition~4 and Lemmas~1--2]{Huster2018Lipschitz}, applied
to that extension. For continuous piecewise affine $g_i$ and real
coefficients $b,a_i$, we have
\begin{align}\label{eq:counting_variation}
    V_*\left(b+\sum_i a_i g_i\right)
        \leq\sum_i|a_i|V_*(g_i),
    \qquad V_*(\sigma\circ g)\leq V_*(g).
\end{align}
The first inequality follows from subadditivity and homogeneity of
variation. For the second, consider the finite sequence of slopes of
the constant extension, including the zero slopes at both ends.
On a component of its nonpositive set with nonempty interior, ReLU
replaces all interior slopes by zero. If the component has positive
regions on both sides, the adjacent slopes satisfy
$s_-\leq0\leq s_+$. The original variation along this block is at
least $|s_+-s_-|=|s_-|+|s_+|$, which is exactly the variation of the
replacement $s_-,0,s_+$. For an unbounded component, use zero for
the absent outer slope; the same argument applies because the extension
has zero derivative there. This also covers components meeting $0$
or $1$. Zero plateaus are included in these blocks. An isolated zero
does not change the derivative almost everywhere, and if the whole
function is nonpositive the new variation is zero. Applying the block
argument to every component proves the second inequality.

Fix the coordinates other than $a$. The varying input $t\mapsto t$
has $V_*=2$, and each constant input has $V_*=0$.
Propagating~\eqref{eq:counting_variation} through the affine layers and
the intervening ReLU activations therefore yields
\begin{align}\label{eq:counting_line_variation}
    V_*\bigl(t\mapsto\realization_\Phi
        (x_1,\ldots,x_{a-1},t,x_{a+1},\ldots,x_d)\bigr)
    \leq2c_a.
\end{align}
This argument also allows the final affine layer to have signed output.

Put $h:=2^{-m}$. If $\tau_\nu\in(0,1)$ are the breakpoints of $g$
and $\Delta_\nu$ the jumps of $g'$, then
\[
    g(t)=g(0)+g'(0+)t+\sum_\nu\Delta_\nu\sigma(t-\tau_\nu),
    \qquad \sum_\nu|\Delta_\nu|=\operatorname{Var}_{(0,1)}(g').
\]
The midpoint rule integrates the affine part exactly. For a single hinge
$t\mapsto\sigma(t-\tau)$ with $\tau\in[jh,(j+1)h]$, it is also exact
on every grid interval other than the one containing $\tau$, since the
hinge is affine there. On the interval containing $\tau$, its error is
\begin{align*}
    \int_{jh}^{(j+1)h}\sigma(t-\tau)\,dt
        -h\sigma((j+1/2)h-\tau) &=\frac{((j+1)h-\tau)^2}{2}
        -h\max\{(j+1/2)h-\tau,0\}\\
    &=\begin{cases}
        (\tau-jh)^2/2,&\tau\leq(j+1/2)h,\\
        ((j+1)h-\tau)^2/2,&\tau\geq(j+1/2)h,
    \end{cases}\\
    &=\frac12\min\{(\tau-jh)^2,((j+1)h-\tau)^2\}
        \leq\frac{h^2}{8}.
\end{align*}
For each breakpoint $\tau_\nu$, choose an index $j_\nu$ such that
$\tau_\nu\in[j_\nu h,(j_\nu+1)h]$. By linearity of integration and
quadrature, followed by the triangle inequality, we therefore obtain
\begin{align*}
    \left|\int_0^1g(t)\,dt
        -h\sum_{j=0}^{2^m-1}g((j+1/2)h)\right| &=\left|\frac12\sum_\nu\Delta_\nu
        \min\{(\tau_\nu-j_\nu h)^2,((j_\nu+1)h-\tau_\nu)^2\}\right|\\
    &\leq\frac12\sum_\nu|\Delta_\nu|
        \min\{(\tau_\nu-j_\nu h)^2,((j_\nu+1)h-\tau_\nu)^2\}\\
    &\leq\frac{h^2}{8}\sum_\nu|\Delta_\nu|
        =\frac{h^2}{8}\operatorname{Var}_{(0,1)}(g').
\end{align*}
Several breakpoints may lie in the same grid interval: each hinge
contributes its own term to this sum. Breakpoints on grid boundaries
contribute zero error.

We now pass from this one-dimensional estimate to the tensor-product
rule. The one-dimensional estimate, together with
\eqref{eq:counting_line_variation}, holds for every fixed choice of
the other coordinates, with the same bound $2c_a$ on the slope
variation. Taking the supremum over those coordinates therefore gives
\[
    \|(I_a-Q_a)\realization_\Phi\|_\infty
    \leq\frac{h^2}{8}\,2c_a=\frac{h^2c_a}{4}.
\]
To combine these bounds, replace integration by quadrature one
coordinate at a time. Subtracting consecutive products
and summing gives the telescoping identity
\begin{align*}
    I_1\cdots I_d-Q_1\cdots Q_d =\sum_{a=1}^d Q_1\cdots Q_{a-1}(I_a-Q_a)I_{a+1}\cdots I_d = \sum_{a=1}^d Q_1\cdots Q_{a-1}I_{a+1}\cdots I_d(I_a-Q_a),
\end{align*}
where the second equality uses the commutativity of the operators. Empty
products are interpreted as the identity operator. Applied to $\realization_\Phi$,
the two products on the left are constant functions equal to the exact
integral and $Q_m(\Phi)$, respectively. The triangle inequality and
the contraction bounds now yield
\begin{align*}
    \left|\int_{[0,1]^d}\realization_\Phi(\bx)\,d\bx-Q_m(\Phi)\right|
    &\leq\sum_{a=1}^d
        \bigl\|Q_1\cdots Q_{a-1}I_{a+1}\cdots I_d(I_a-Q_a)\realization_\Phi
            \bigr\|_\infty\leq\sum_{a=1}^d\|(I_a-Q_a)\realization_\Phi\|_\infty\\
    &\leq\frac{h^2}{4}\sum_{a=1}^dc_a
        \leq\frac{K_\Phi}{4}\,2^{-2m},
\end{align*}
which is~\eqref{eq:counting_midpoint_error}. The argument uses only
the slope jumps of the piecewise affine coordinate restrictions, not
second-order differentiability.
\end{proof}

For the counting reduction, let $\Phi$ be
scalar-valued, with actual depth $L$, and let
$\theta\subset\Q$ be a finite set of parameters, i.e., each entry of $\bW_\ell$ and $\bb_\ell$ is contained in $\theta$. Write each parameter in reduced form
as $a/b$, with $a\in\Z$ and $b\in\N$, normalizing zero to $0/1$.
Choose $p\in\N$ such that $|a|\leq2^p$ and $b\leq2^p$ for every
parameter. For fixed $\theta$, this choice is independent of the
architecture and quadrature accuracy. Define
\begin{align*}
    P_\ell&:=\#\operatorname{supp}(\bW_\ell)+w_\ell,
    &r_\ell&:=\max\left\{1,\max_i
        \#\{j:(\bW_\ell)_{ij}\neq0\}\right\},\\
    P_\Phi&:=\sum_{\ell=1}^{L}P_\ell,
    &r_\Phi&:=\max_{\ell=1,\ldots,L}r_\ell.
\end{align*}
Thus $P_\Phi$ counts nonzero weights and all bias positions, and
$r_\Phi$ bounds the nonzero inputs into a layer.

\begin{theorem}\label{thm:counting_quadrature}
Let $\Phi$ be a scalar-valued ReLU network of depth $L$ whose
parameters are rational numbers $a/b$, with $|a|\leq2^p$ and
$1\leq b\leq2^p$, as above. The denominators may differ between
parameters, including within a layer. For each layer $\ell$, let
$q_\ell$ be the least common multiple of its parameter denominators.
Then $q_\ell\bW_\ell$ and $q_\ell\bb_\ell$ have integer entries.
Let $p_\ell$ be the smallest positive integer for which $q_\ell$ and
the absolute values of all these integer entries are at most
$2^{p_\ell}$. For $m\in\N$, put
\begin{align}\label{eq:counting_size}
    B&:=m+3+\sum_{\ell=1}^{L}
        \bigl(p_\ell+\lceil\log_2(r_\ell+1)\rceil\bigr),\nonumber\\
    S_{\Phi,m}&:=dm+B\sum_{\ell=1}^{L}P_\ell p_\ell.
\end{align}
There is a uniform algorithm which, given $\Phi,m,p$, constructs a
$3$-CNF formula $\alpha_m$ whose number of satisfying assignments
determines the midpoint quadrature sum from~\eqref{eq:counting_midpoints}
exactly:
\begin{align}\label{eq:counting_signed_sum}
    Q_m(\Phi)=
    \frac{\#\mathrm{SAT}(\alpha_m)-2^{md+B-1}}
         {2^{(d+1)m+1}\prod_{\ell=1}^{L}q_\ell}.
\end{align}
The formula has $n_m=\mathcal O(S_{\Phi,m})$ variables, at most $4n_m$
clauses, and binary encoding length
$\mathcal O(S_{\Phi,m}\log_2(2+S_{\Phi,m}))$.
In terms of the original parameter bound $p$, we have
\[
    S_{\Phi,m}=\mathcal O\bigl(dm+P_\Phi^2p(m+P_\Phi p)\bigr).
\]

Constructing $\alpha_m$ takes
$\operatorname{poly}(B_{\rm in},S_{\Phi,m})$ bit operations, where
$B_{\rm in}$ is the total binary encoding length of $\Phi,m,p$.
Once the exact integer $\#\mathrm{SAT}(\alpha_m)$ is supplied,
recovering the rational number $Q_m(\Phi)$
by~\eqref{eq:counting_signed_sum} takes polynomially many bit operations
in the same quantities. These bounds do not include the cost of
computing $\#\mathrm{SAT}(\alpha_m)$. The resulting quadrature value
$Q_m(\Phi)$ satisfies the integration error bound
\eqref{eq:counting_midpoint_error}.

\end{theorem}
\begin{proof}
We first construct a Boolean circuit which computes an exact integer
numerator of $\realization_\Phi(\bx_{\boldsymbol j})$ from the binary
digits of a midpoint index $\boldsymbol j\in J_{m,d}$. We then add
input bits and a comparison test so that counting the accepted inputs
recovers the sum of these numerators. Finally, we encode this Boolean
circuit by a $3$-CNF formula, preserving the number of accepted inputs.

\emph{Exact integer evaluation of the network.}
Reduce the supplied fractions and normalize zero to $0/1$, if
necessary. This cannot increase their numerator or denominator bounds.
Then compute $q_\ell$ and multiply every parameter of layer $\ell$
by $q_\ell$. These integer coefficients represent the same layer
when accompanied by the denominator $q_\ell$. The least common
multiples, exact quotients, and scaled coefficients are computed
before building the Boolean circuit. The bound on their bit lengths
will be verified below.
Use $m$ free input bits, each allowed to take either value $0$ or $1$
independently, for each coordinate index
$j_a\in\{0,\ldots,2^m-1\}$, hence $dm$ grid bits in total.
Each assignment of these bits specifies exactly one midpoint.
Appending a fixed binary one to the bits of $j_a$ gives the numerator
$2j_a+1$ of $x_{\boldsymbol j,a}$ with denominator $D_0:=2^{m+1}$.
At layer $\ell$, use the common positive denominator
$D_\ell:=2^{m+1}\prod_{s=1}^\ell q_s$, so that
$D_\ell=q_\ell D_{\ell-1}$. If the inputs to a neuron are
$u_s/D_{\ell-1}$, its nonzero weights, after this rewriting, are
$a_s/q_\ell$, and its bias is $c/q_\ell$, then its preactivation is
\[
    \sum_s\frac{a_s}{q_\ell}\frac{u_s}{D_{\ell-1}}
        +\frac c{q_\ell}
    =\frac{\sum_s a_su_s+cD_{\ell-1}}{D_\ell}.
\]
Thus the preactivation numerator is the integer
\begin{align}\label{eq:counting_integer_layer}
    v=\sum_s a_su_s+cD_{\ell-1}.
\end{align}
At a hidden neuron, the activation numerator is $\max\{v,0\}$,
since $D_\ell>0$. At the final affine layer, retain $v$ without
applying ReLU. Induction over the layers therefore computes an integer
$U_{\boldsymbol j}$ satisfying
\[
    \realization_\Phi(\bx_{\boldsymbol j})
        =\frac{U_{\boldsymbol j}}{D_{L}}
        =\frac{U_{\boldsymbol j}}
            {2^{m+1}\prod_{\ell=1}^{L}q_\ell}.
\]
No rounding or rational division is needed in this evaluation.

\emph{Bit width and Boolean circuit size.}
The bit width $B$ is the number of binary digits used to represent
each intermediate integer. To specify how negative integers are
represented, let $\beta_{B-1},\ldots,\beta_0\in\{0,1\}$ be these
digits, written from most to least significant. Their usual
\emph{unsigned} interpretation is $\sum_{k=0}^{B-1}2^k\beta_k$.
Their \emph{two's-complement} interpretation instead assigns a negative
weight to the first digit:
\[
    -2^{B-1}\beta_{B-1}+\sum_{k=0}^{B-2}2^k\beta_k.
\]
This represents each integer in $\{-2^{B-1},\ldots,2^{B-1}-1\}$
uniquely. The first digit is called the \emph{sign bit}: it is $1$
exactly for negative integers. In general, the unsigned interpretation of the representation of $U$
is $U$ for $U\geq0$ and $2^B+U$ for $U<0$.

Set $M_0:=2^{m+1}$ and
$M_\ell:=2^{m+1}\prod_{s=1}^\ell[2^{p_s}(r_s+1)]$.
The input numerators have magnitude less than $M_0$, and
$D_\ell\leq M_\ell$. Inductively, the sum of the absolute magnitudes
of the terms in~\eqref{eq:counting_integer_layer} is at most
\[
    r_\ell2^{p_\ell}M_{\ell-1}
        +2^{p_\ell}D_{\ell-1}\leq M_\ell.
\]
Thus, every activation numerator, product, and affine partial sum has
magnitude at most $M_\ell\leq M_{L}\leq2^{B-2}$.
For multiplication by a fixed coefficient, first multiply absolute
values by shifts and additions, then restore the sign. All partial
sums in this multiplication are nonnegative and bounded by the absolute
product, so the same $B$ signed bits suffice for these intermediate
operations as well. All these values and their negatives lie in the
representable range. Thus, there is no overflow, meaning that no exact
result lies outside the range of the chosen representation.
Preparing the constants $cD_{\ell-1}$ uses polynomial bit work.

In this representation, ordinary binary addition, retaining the $B$
least significant bits, also gives signed addition whenever the exact
result is in the representable range. At each digit, the sum digit and
the carry to the next position depend on only two input digits and the
previous carry. A fixed number of Boolean gates suffices per position,
so addition uses $\mathcal O(B)$ gates. Negation replaces every bit
$\beta_k$ by $1-\beta_k$ and adds one, again retaining $B$ bits, and
has the same gate bound.
A sign test reads the sign bit. The ReLU operation keeps all bits when
that bit is zero and sets them all to zero otherwise. It too uses
$\mathcal O(B)$ gates.
For multiplication by a fixed coefficient, expand its absolute value
in binary. Each nonzero digit requires one copy of the absolute input
shifted by the corresponding number of positions. Such a shift simply
multiplies its value by a power of two. Adding these at most
$\mathcal O(p_\ell)$ copies in layer $\ell$ and restoring the sign
therefore uses $\mathcal O(p_\ell B)$ gates of bounded fan-in.
Only nonzero connections need multipliers. Including biases and all
neurons gives $\mathcal O(B\sum_{\ell=1}^{L}P_\ell p_\ell)$ gates
for exact forward evaluation.
The circuit outputs the two's-complement representation of
$U_{\boldsymbol j}$ in $B$ bits. The magnitude bound above ensures
\[
    -2^{B-1}\leq U_{\boldsymbol j}<2^{B-1}.
\]

\emph{Encoding the numerator sum as a count.}
To turn these possibly signed numerators into nonnegative integers
that can be counted by a threshold test, flip the sign bit and then
interpret the resulting bits as unsigned. If $U_{\boldsymbol j}\geq0$,
the sign bit changes from $0$ to $1$, adding $2^{B-1}$ to the unsigned
value $U_{\boldsymbol j}$. If $U_{\boldsymbol j}<0$, the sign bit changes
from $1$ to $0$, subtracting $2^{B-1}$ from the unsigned value
$2^B+U_{\boldsymbol j}$. In both cases the resulting unsigned integer is
\[
    \widetilde U_{\boldsymbol j}:=U_{\boldsymbol j}+2^{B-1}
        \in\{0,\ldots,2^B-1\}.
\]
Extend the evaluation circuit by $B$ additional free input bits
encoding an unsigned threshold $t\in\{0,\ldots,2^B-1\}$ and by a
\emph{comparator}, a Boolean circuit whose output is one exactly when
$t<\widetilde U_{\boldsymbol j}$. To compare the two unsigned integers,
inspect their digits from most to least significant: the first
differing digit decides which is smaller. Maintaining whether the
digits seen so far are equal and, if not, which integer is smaller
requires only a fixed number of gates per position. Thus the comparator
uses $\mathcal O(B)$ gates. Its output is the single output bit of
the extended circuit. For fixed $\boldsymbol j$, the
accepted thresholds are precisely
$0,\ldots,\widetilde U_{\boldsymbol j}-1$, with none when
$\widetilde U_{\boldsymbol j}=0$. Consequently, the number of accepted
assignments of all $dm+B$ free input bits is
$\sum_{\boldsymbol j\in J_{m,d}}\widetilde U_{\boldsymbol j}$.

\emph{Conversion to a counting-preserving CNF formula.}
Replace the gates of the extended circuit by AND and NOT gates,
incurring only a constant-factor increase in size. Introduce one
Boolean variable for each gate output and each constant wire, in
addition to the $dm+B$ free input variables. Impose full Tseytin
equivalences~\cite{Tseitin1983,sipser13} to enforce the value computed
by each gate. Constant wires are fixed by
\emph{unit clauses}, each consisting of just one literal: $(z)$ forces
$z=1$, and $(\neg z)$ forces $z=0$. Add a unit clause requiring the
comparator output to be one, and call the resulting conjunction
$\alpha_m$.
For each assignment of the free input variables, the constants and
then the gates, evaluated in topological order (each gate after all of
its inputs), determine all auxiliary variables uniquely. This unique
extension satisfies the gate and
constant clauses, and it satisfies the output clause if and only if
$t<\widetilde U_{\boldsymbol j}$. Thus each accepted pair
$(\boldsymbol j,t)$ gives exactly one satisfying assignment of
$\alpha_m$, and every satisfying assignment arises in this way.

To ensure that all free inputs occur in the formula, also add the
tautological clause $(x\lor\neg x)$ for each grid and threshold bit.
These clauses impose no restrictions and do not introduce variables.
They ensure that even unused input coordinates or a constant network
do not cause free bits to disappear from the count. It follows that
\[
    \#\mathrm{SAT}(\alpha_m)
        =\sum_{\boldsymbol j\in J_{m,d}}\widetilde U_{\boldsymbol j}
        =\sum_{\boldsymbol j\in J_{m,d}}U_{\boldsymbol j}+2^{md+B-1}.
\]
Since $Q_m(\Phi)=2^{-md}D_{L}^{-1}
\sum_{\boldsymbol j\in J_{m,d}}U_{\boldsymbol j}$ and
$2^{md}D_{L}=2^{(d+1)m+1}\prod_{\ell=1}^{L}q_\ell$, this proves
\eqref{eq:counting_signed_sum}.

\emph{Formula size, construction, and recovery of the quadrature sum.}
There are $dm+B$ free variables and
$\mathcal O(B\sum_{\ell=1}^{L}P_\ell p_\ell+B)$ gate and constant
variables. By~\eqref{eq:counting_size}, their total number $n_m$ is
$\mathcal O(S_{\Phi,m})$.
Every gate contributes at most three clauses, each constant one clause,
and each free variable one tautology. Including the output clause,
the number of clauses is at most four times the actual number of
variables. Each clause has at most three literals. Variable labels use
$\mathcal O(\log_2(2+S_{\Phi,m}))$ bits, proving the encoding bound.
For each layer, the least common multiple $q_\ell$ is at most the
product of the $P_\ell$ parameter denominators (zero weights need not
be included). Thus, $q_\ell\leq2^{P_\ell p}$. If a parameter in that
layer is $a/b$, then $q_\ell/b$ is at most the product of the other
denominators, so $|a|q_\ell/b\leq2^{P_\ell p}$. Consequently,
\[
    p_\ell\leq P_\ell p,\qquad
    \sum_{\ell=1}^{L}P_\ell p_\ell\leq P_\Phi^2p,
    \qquad B\leq m+3+P_\Phi p+P_\Phi.
\]
For the last inequality we used
$\lceil\log_2(r_\ell+1)\rceil\leq r_\ell\leq P_\ell$.
These bounds prove the asserted estimate for $S_{\Phi,m}$ in terms
of the original parameter bit bound $p$.
The layer denominators can be formed by repeated use of
$\operatorname{lcm}(u,v)=uv/\gcd(u,v)$, with the greatest common
divisors computed by the Euclidean algorithm. These operations and the
subsequent exact divisions use polynomially many bit operations on
integers of the bounded lengths above.
All gates and clauses are generated by explicit loops over neurons,
stored connections, and bit positions. The $2^{md}$ midpoint indices
are represented by free bits, not enumerated. Reading the network,
computing the layer denominators and scaled coefficients, forming
$D_\ell$ and $cD_{\ell-1}$, and writing the
formula therefore take
$\operatorname{poly}(B_{\rm in},S_{\Phi,m})$ bit operations.

After the integer $\#\mathrm{SAT}(\alpha_m)$ has been computed by a
separate counting algorithm, recover $Q_m(\Phi)$ by the subtraction
and fraction in~\eqref{eq:counting_signed_sum}. The count has at most
$n_m+1=\mathcal O(S_{\Phi,m})$ bits. The correction integer
$2^{md+B-1}$ and the denominator
$2^{(d+1)m+1}\prod_{\ell=1}^{L}q_\ell$ also have
$\mathcal O(S_{\Phi,m})$ bits. Thus this recovery takes polynomially
many bit operations and may return an exact numerator--denominator
pair without reducing the fraction. It introduces no additional error and
the integration error is precisely the midpoint quadrature error
bounded in~\Cref{lem:counting_midpoint_error}.
\end{proof}

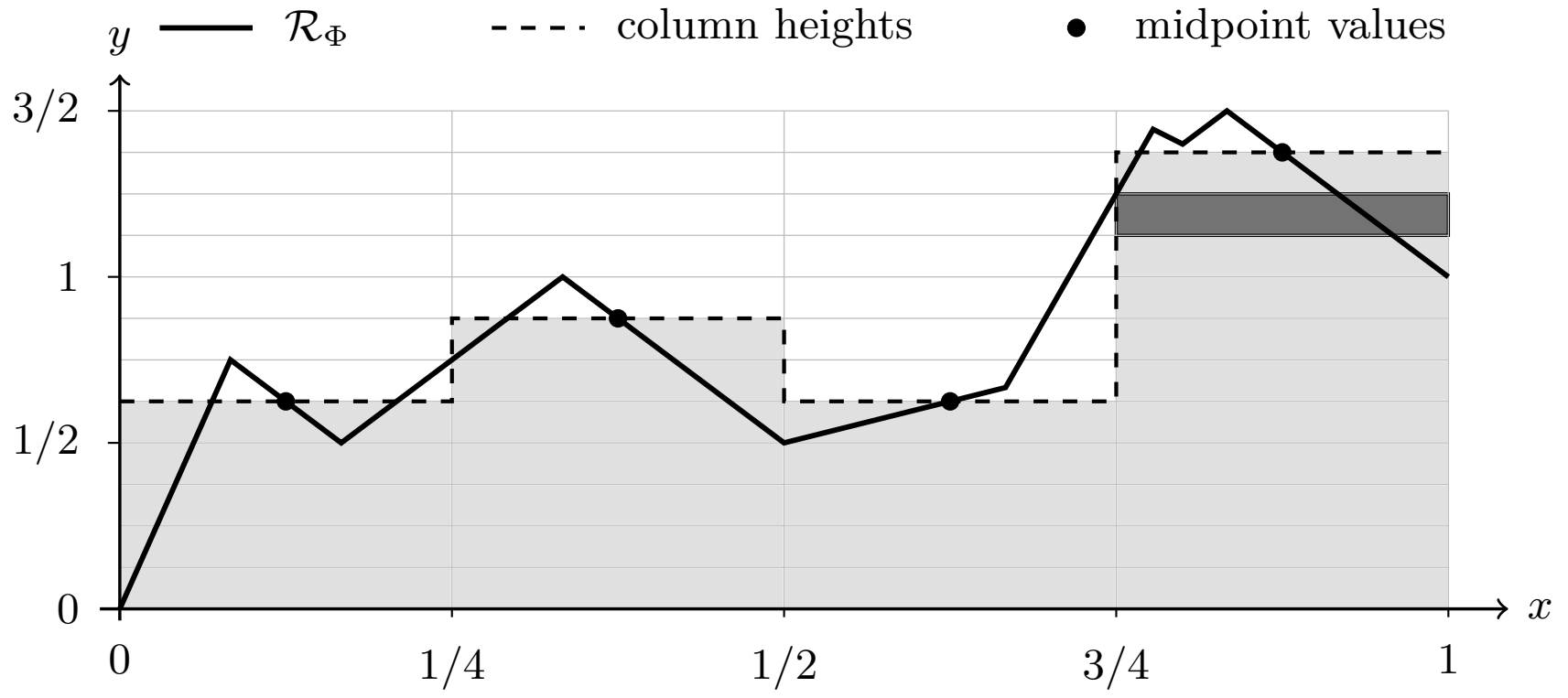
\begin{figure}[htbp]
    \centering
    \input{counting_rectangles.tex}
    \caption{Rectangle counting for
    $\realization_\Phi(x)=x+g(g(v(x)))$ on $[0,1]$, where
    $v(x):=\sigma(2x)-\sigma(5x-3)+\sigma(5x-4)$ and
    $g(t):=\sigma(2t)-\sigma(3t-1)$.
    Here $m=2$ and $q_\ell=1$ for every layer.}
    \label{fig:counting_rectangles}
\end{figure}

\begin{remark}[Nonnegative realizations]
If $\realization_\Phi\geq0$ on $[0,1]^d$ is known in advance, the
construction in~\Cref{thm:counting_quadrature} simplifies. Every output
numerator $U_{\boldsymbol j}$ is then nonnegative. Omit the sign-bit
flip and let the comparator test $t<U_{\boldsymbol j}$ instead of
$t<U_{\boldsymbol j}+2^{B-1}$. There are exactly
$U_{\boldsymbol j}$ accepted thresholds for each midpoint, and the
same counting-preserving CNF encoding therefore gives
\[
    \#\mathrm{SAT}(\alpha_m)
        =\sum_{\boldsymbol j\in J_{m,d}}U_{\boldsymbol j},
    \qquad
    Q_m(\Phi)=\frac{\#\mathrm{SAT}(\alpha_m)}
        {2^{(d+1)m+1}\prod_{\ell=1}^{L}q_\ell}.
\]
The size and construction bounds are unchanged, and no correction term
is needed. This variant uses the supplied nonnegativity information and
it does not require the algorithm to decide whether the network is
nonnegative. Figure~\ref{fig:counting_rectangles} illustrates this case
as an area made up of exact counting rectangles, whose column heights
agree with the midpoint values without output rounding.
\end{remark}

The sparse estimate counts only stored nonzero connections and all
neurons. If the input is dense, reading its zero entries is included
separately through $B_{\rm in}$. A sharper bound applies when all
parameters can be written with one common denominator $q$, where
$q\leq2^p$ and all corresponding numerators have absolute value at
most $2^p$. Each $q_\ell$ then divides $q$, so the scaled coefficients
also satisfy $p_\ell\leq p$. Consequently,
\[
    S_{\Phi,m}\leq dm+P_\Phi p
        \left[m+3+L\bigl(p+\lceil\log_2(r_\Phi+1)\rceil\bigr)\right].
\]
In this common-denominator case, $d\leq w$,
$P_\Phi\leq L(w^2+w)$ and $r_\Phi\leq w$, giving the simpler dense bound
\[
    S_{\Phi,m}=\mathcal O\bigl(dm+Lw^2p[m+1+L(p+\log_2(w+1))]\bigr).
\]

\begin{corollary}\label{cor:counting_complexity}
Assume the deterministic exact counting bound of
\cite[Theorem~1.1]{OpenAI2026ETH}.
There is an absolute integer $C_0\geq1$ such that, for every finite
$\theta\subset\Q$, there is a uniform deterministic
algorithm $\cA_{\#\mathrm{SAT}}$ with the following property.
Let $q$ be the least common multiple of the reduced denominators in
$\theta$, and fix $p\in\N$ such that $q\leq2^p$ and
$|q\vartheta|\leq2^p$ for every $\vartheta\in\theta$.
For every scalar-valued ReLU network $\Phi$ with parameters in $\theta$,
put
\begin{align}\label{eq:counting_work_constant}
    c_\Phi:=\frac{1}{16C_0(1+d+P_\Phi p)}.
\end{align}
There is $N_{0,\Phi}\in\N$, with
$N_{0,\Phi}\geq16$, such that, for every integer $N\geq N_{0,\Phi}$,
the algorithm returns a rational $\cA_{\#\mathrm{SAT}}(\Phi,N)$ in
total word-RAM runtime at most $N$ and satisfies
\begin{align}\label{eq:counting_work_error}
    \left|\cA_{\#\mathrm{SAT}}(\Phi,N)
        -\int_{[0,1]^d}\realization_\Phi(x)\,dx\right|
    \leq N^{-c_\Phi\sqrt{\log_2\log_2N}}.
\end{align}
Here $C_0$ is a fixed constant in the complete counting runtime bound,
independent of the network and the work budget. The threshold
$N_{0,\Phi}$ may depend on the network encoding. 
\end{corollary}
\begin{proof}
The cited theorem states that for fixed $k$ and $D$, exact counting on
$k$-CNFs with $n\geq20$ declared variables and $M\leq Dn$ clauses
uses $(n+M)^{\mathcal O_{k,D}(1)}2^{\mathcal O_{k,D}(n/\sqrt{\log n})}$
bit operations and polynomial space. Apply it with $k=3$ and $D=4$
to \Cref{thm:counting_quadrature}. Since $n=\mathcal O(S_{\Phi,m})$, the complete cost of computing
$Q_m(\Phi)$ is at most
\begin{align}\label{eq:counting_runtime}
    \operatorname{poly}(B_{\rm in},S_{\Phi,m})
    \exp\left(\mathcal O\left(
       \frac{S_{\Phi,m}}{\sqrt{\log_2(2+S_{\Phi,m})}}\right)\right),
\end{align}
including input, construction, exact counting, and rational output.
Here $B_{\rm in}$ includes the binary encoding of the budget $N$.

We give an integer choice of precision to invert this bound.
The common-denominator choice of $p$ in the statement gives
$p_\ell\leq p$ for every layer. The size parameter therefore satisfies
\[
    S_{\Phi,m}\leq(d+P_\Phi p)m+
    P_\Phi p[3+L(p+\lceil\log_2(r_\Phi+1)\rceil)].
\]
Choose the absolute integer $C_0\geq1$ so that the base-two logarithm of the
complete runtime bound is at most
\[
    C_0\left[\frac{S_{\Phi,m}}{\sqrt{\log_2(2+S_{\Phi,m})}}
        +\log_2(2+B_{\rm in}+S_{\Phi,m})\right].
\]
Increasing $C_0$ also accounts for computing the precision and the
finite-word simulation described below. This constant is fixed for
the counting algorithm, independently of $\Phi$ and $N$.
With $c_\Phi$ from~\eqref{eq:counting_work_constant}, for $N\geq16$, set
\begin{align}\label{eq:counting_budget_precision}
    b_N&:=\lfloor\log_2N\rfloor,
    &t_N&:=\left\lfloor\sqrt{\lfloor\log_2b_N\rfloor}\right\rfloor,
    &m_N&:=\max\{1,\lfloor c_\Phi b_Nt_N\rfloor\}.
\end{align}
Bit lengths, integer square roots, and integer division compute these
quantities exactly and no real logarithm or square root is evaluated.
The algorithm returns $Q_{m_N}(\Phi)$.

For fixed $\Phi$, as $N\to\infty$, we have
$m_N\sim c_\Phi\log_2N\sqrt{\log_2\log_2N}$,
$\log_2(2+S_{\Phi,m_N})\sim\log_2\log_2N$, and
$B_{\rm in}=\mathcal O_\Phi(\log N)$.
Consequently, the logarithm of the total runtime is at most
$C_0(d+P_\Phi p)c_\Phi\log_2N+o_\Phi(\log N)$, hence at most
$\frac14\log_2N$ for all sufficiently large $N$.
In particular, the complete word-RAM runtime is at most $N$.
Finally, \Cref{lem:counting_midpoint_error} gives an error bounded by
$K_\Phi2^{-2m_N}/4$. The asymptotic formula for $m_N$ absorbs this
fixed prefactor and yields~\eqref{eq:counting_work_error}, after
increasing $N_{0,\Phi}$.
If $K_\Phi=0$, the realization is constant and one exact evaluation
suffices. On a fixed class $\nn(d,w,L,\theta)$, both the encoding length
and $d+P_\Phi p$ are uniformly bounded, and
$K_\Phi\leq(r_\Phi2^p)^{L}$ is uniformly bounded as well.
The constants can therefore be chosen uniformly on that class.
\end{proof}

\section{Numerical experiments}\label{sec:rates}
\input{benchmark_results.tex}

\printbibliography
\end{document}

%% file: mymacros.tex
\def\R{{\mathbb R}}

\def\N{{\mathbb N}}
\def\Z{{\mathbb Z}}

\newcommand{\cA}{\mathcal{A}}

\def\bx{{\boldsymbol{x}}}

\def\bb{{\boldsymbol{b}}}

\def\bW{{\boldsymbol{W}}}
\def\bb{{\boldsymbol{b}}}

\def\norm#1#2{\|#1\|_{#2}}

\def\v{\mathbf{v}}

\def\normL2#1#2{\|#1\|_{L^2(#2)}}

\newtheorem{theorem}{Theorem}[section]

\newtheorem{lemma}[theorem]{Lemma}
\newtheorem{corollary}[theorem]{Corollary}

\newenvironment{remark}[1][]{%
    \medskip
    \noindent
    \textbf{Remark}%
    \if\relax\detokenize{#1}\relax
    \else\ \textnormal{(#1)}%
    \fi
    \textbf{.}\itshape
}{%
    \qed
    \smallskip
}

\DeclareRobustCommand{\CSrange}[1]%
{\ensuremath{[#1]}}

%% file: counting_rectangles.tex
\begin{tikzpicture}[x=9.6cm,y=2.4cm,
    every node/.style={font=\footnotesize},
    every path/.style={draw=black}]
    \foreach \j/\height in {0/5,1/7,2/5,3/11}{%
        \pgfmathtruncatemacro{\lastlevel}{\height-1}
        \foreach \t in {0,...,\lastlevel}{%
            \path[fill=black!12,draw=white,line width=0.35pt]
                ({\j/4},{\t/8}) rectangle ({(\j+1)/4},{(\t+1)/8});
        }
    }
    \path[fill=black!55,draw=black!80!black,line width=0.65pt]
        ({3/4},{9/8}) rectangle (1,{10/8});
    \foreach \j in {0,...,4}
        \draw[black!25,line width=0.25pt] ({\j/4},0) -- ({\j/4},1.5);
    \foreach \t in {0,...,12}
        \draw[black!25,line width=0.25pt] (0,{\t/8}) -- (1,{\t/8});
    \draw[black!65!black,dashed,line width=0.8pt]
        (0,{5/8})
        \foreach \j/\height in {0/5,1/7,2/5,3/11}{%
            -- ({\j/4},{\height/8}) -- ({(\j+1)/4},{\height/8})
        };
    \draw[black,line width=1.1pt]
        (0,0) -- ({1/12},{3/4}) -- ({1/6},{1/2}) -- ({1/3},1)
        -- ({1/2},{1/2}) -- ({2/3},{2/3}) -- ({7/9},{13/9})
        -- ({4/5},{7/5}) -- ({5/6},{3/2}) -- (1,1);
    \foreach \j/\height in {0/5,1/7,2/5,3/11}
        \fill[black] ({(\j+0.5)/4},{\height/8}) circle[radius=1.7pt];
    \draw[->,line width=0.65pt] (-0.015,0) -- (1.045,0)
        node[right] {$x$};
    \draw[->,line width=0.65pt] (0,-0.035) -- (0,1.61)
        node[above] {$y$};
    \foreach \v/\lab in {0/0,0.25/{1/4},0.5/{1/2},0.75/{3/4},1/1}{%
        \draw[line width=0.4pt] (\v,0) -- (\v,-0.025)
            node[below=2pt] {$\lab$};
    }
    \foreach \v/\lab in {0/0,0.5/{1/2},1/1,1.5/{3/2}}{%
        \draw[line width=0.4pt] (0,\v) -- (-0.009,\v)
            node[left=2pt] {$\lab$};
    }
    \draw[line width=1.1pt] (0.03,1.75) -- (0.10,1.75);
    \node[anchor=west] at (0.11,1.75) {$\realization_\Phi$};
    \draw[black!65!black,dashed,line width=0.8pt] (0.28,1.75) -- (0.35,1.75);
    \node[anchor=west] at (0.36,1.75) {column heights};
    \fill[black] (0.72,1.75) circle[radius=1.7pt];
    \node[anchor=west] at (0.75,1.75) {midpoint values};
\end{tikzpicture}

%% file: benchmark_results.tex
All numerical examples use eight or twelve input dimensions. SAT denotes the model counting quadrature algorithm developed above,  MC denotes uniform Monte
Carlo. Accuracy is reported only as relative error. For a reference $I>0$,
the relative error is $|Q-I|/I$; for repeated MC trials we report
$\sqrt{r^{-1}\sum_i(Q_i-I)^2}/I$.

The implementation is written in Python and uses the GANAK model
counter~\cite{Sharma2019Ganak,Soos2025Ganak} through the Python package
\texttt{pyganak}, version~2.9.0. Thus the
experiments use GANAK's practical counting implementation rather than the
deterministic algorithm assumed in Corollary~\ref{cor:counting_complexity}.
The measurements were obtained on an Apple M3 with 24\,GB of memory under
macOS~14.8.9, using Python~3.12.7 and NumPy~2.5.5, with numerical-library
thread counts restricted to one.

The benchmark implementation uses a left-endpoint tensor grid and
floor quantization of the output. For $F=\realization_\Phi$ with
$0\leq F<1$, the reported estimate is
\[
    \widehat I_{m,k}
    =2^{-dm-k}\sum_{\boldsymbol j\in\{0,\ldots,2^m-1\}^d}
        \left\lfloor 2^kF(2^{-m}\boldsymbol j)\right\rfloor.
\]
Here $k$ is the output precision listed in Table~\ref{tab:bumps}; the
experiments use $m=k+4$ for the corner bump in eight dimensions and
$m=k+5$ for the other two networks. The integer sum is obtained from a
single $3$-CNF model count, using rational network arithmetic and a
$k$-bit threshold comparison, without enumerating the grid points.
This discretization differs from the midpoint rule analyzed above.

Reported times are median wall times over three SAT runs and eleven
independent MC runs per budget, after warm-up. SAT timings include network
normalization, CNF construction, solver initialization, clause transfer,
model counting, and rational postprocessing. 

We test integrands that are supported on a tiny subdomain (bumps) and thus are difficult for MC. The corner bumps are $\sigma(1/2-\sum_{a=1}^{d}x_a)$ in $d=8,12$.
The centered 8D bump is $\sigma(1/4-\sum_{a=1}^{8}|x_a-1/2|)$.
Every coordinate affects the output. The corner network has one hidden
neuron; the centered network has hidden widths 16 and 1.

Table~\ref{tab:bumps} pairs the highest tested SAT precision with the MC
budget closest in median wall time on a logarithmic scale. These are
actual measured settings, without interpolation. The figure shows every
budget; shading is the interquartile range of MC relative errors. All MC trials on the 12D corner returned
zero, hence relative error one. 
\begin{table}[htbp]
\centering\small
\caption{High-dimensional bump comparison. Errors are relative RMSE in
percent and $N=2^p$ for MC.}
\label{tab:bumps}
\begin{tabular}{lrrrrrr}
\toprule
& \multicolumn{3}{c}{SAT} & \multicolumn{3}{c}{Uniform MC}\\
\cmidrule(lr){2-4}\cmidrule(lr){5-7}
Network & $k$ & seconds & error (\%) & $p$ & seconds & error (\%)\\
\midrule
Corner bump 8D & 8 & 3.87 & 1.52 & 27 & 5.13 & 47.6\\
Corner bump 12D & 5 & 6.11 & 23.7 & 27 & 7.52 & 100\\
Centered bump 8D & 5 & 5.41 & 46 & 27 & 6.01 & 24.3\\
\bottomrule
\end{tabular}
\end{table}
\begin{figure}[htbp]
\centering
\includegraphics[width=\linewidth,trim={0cm 0 0 2cm},clip]{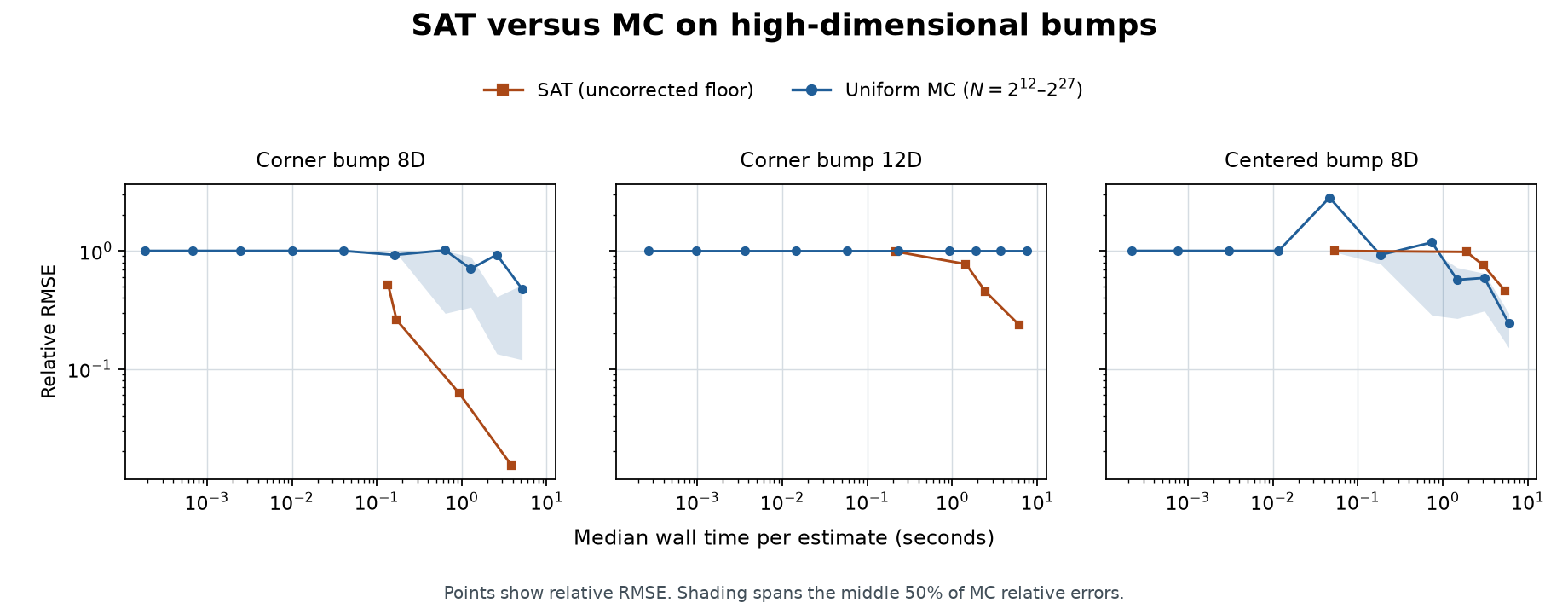}
\caption{Relative RMSE versus median wall time for SAT (\textcolor{red}{red line}) and
uniform MC (\textcolor{blue}{blue line}). }
\end{figure}
\clearpage